\documentclass[5p]{elsarticle}
\usepackage{geometry}
\usepackage{lineno,hyperref}
\usepackage{amsmath}
\usepackage{amsfonts} 
\usepackage{algorithm}
\usepackage[noend]{algpseudocode}
\usepackage{balance}

\makeatletter
\def\ps@pprintTitle{%
  \let\@oddhead\@empty
  \let\@evenhead\@empty
  \let\@oddfoot\@empty
  \let\@evenfoot\@oddfoot
}
\makeatother

\begin{document}

\begin{frontmatter}

\title{Casting graph isomorphism as a point set registration problem \\ using a simplex embedding and sampling}

\author{Yigit Oktar}

\address{orcid: 0000-0002-8736-8013}

\begin{abstract}
Graph isomorphism is an important problem as its worst-case time complexity is not yet fully understood. In this study, we try to draw parallels between a related optimization problem called point set registration. A graph can be represented as a point set in enough dimensions using a simplex embedding and sampling. Given two graphs, the isomorphism of them corresponds to the existence of a perfect registration between the point set forms of the graphs. In the case of non-isomorphism, the point set form optimization result can be used as a distance measure between two graphs having the same number of vertices and edges. The related idea of equivalence classes suggests that graph canonization may be an important tool in tackling graph isomorphism problem and an orthogonal transformation invariant feature extraction based on this high dimensional point set representation may be fruitful. The concepts presented can also be extended to automorphism, and subgraph isomorphism problems and can also be applied on hypergraphs with certain modifications.
\end{abstract}


\end{frontmatter}

\section{Introduction}

Graph isomorphism (GI), a concept in graph theory, questions whether two graphs $G_{1}$ and $G_{2}$ can be related with a bijection through their vertices. Namely, the graphs are isomorphic if and only if one can put the vertices of two graphs into a one-to-one correspondence preserving the edge information. GI requires special attention as its worst-case computational complexity is not yet fully understood. Namely, GI is in class nondeterministic polynomial time (NP). In other words, any candidate solution can be verified in polynomial time by a deterministic Turing machine. However, it is not yet known whether GI is NP-complete or not, and its answer has important consequences~\cite{SCHONING1988312}. It is recently shown that graph isomorphism problem is solvable in quasipolynomial time meaning that "GI is not NP-complete unless all of NP can be solved in quasipolynomial time"~\cite{1512.03547}

The aim of this paper is not to provide an explicit algorithm for GI, but to draw parallels between a related problem, namely point set registration.
Point set registration is the optimization process of aligning two point clouds through finding a transformation between, performed usually in 2D or 3D settings. However, the problem can be thought in arbitrary dimensions. In case of unknown point correspondences, registration of both pose (transformation) and correspondences is needed. Throughout this paper, only rigid transformations are considered. However, non-rigid registration of points is also possible~\cite{5432191}.

\begin{equation}
\label{op}
\underset{\mathbf{M}\in\mathbb{R}^{d\times d}}{\arg\min}\lVert \mathbf{MX}-\mathbf{Y} \rVert_{2}^{2} \quad s.t. \quad \mathbf{M}^{T}\mathbf{M} = \mathbf{I}
\end{equation}

A special case of point set registration is the orthogonal Procrustes optimization problem formulated as in Eqn. (\ref{op}), where $\mathbf{X} \in \mathbb{R}^{d\times n}$ and $\mathbf{Y} \in \mathbb{R}^{d\times n}$ are two sets of $d-dimensional$ points with known correspondences. The transformation matrix is only allowed to be an orthogonal matrix. In that case, $d\times d$ orthogonal matrices form the orthogonal group $O(d)$. If the determinant of $\mathbf{M}$ is also forced to $+1$, the subgroup called the special orthogonal group (or the rotation group) $SO(d)$ is attained, generalizing the rotation transformation to $d$ dimensions. All other matrices have determinant $-1$. A closed form solution to the orthogonal Procrustes problem is possible through singular value decomposition~\cite{Schonemann1966}.

On the dual side, if the transformation is known and the correspondences are needed the formulation then is as in Eqn. (\ref{perm}), where $\mathbf{P}$ is a permutation matrix from the set of all $n\times n$ permutation matrices $\mathbb{\mathcal{P}}_{n}$.

\begin{equation}
\label{perm}
\underset{\mathbf{P}\in\mathbb{\mathcal{P}}_{n}}{\arg\min}\lVert \mathbf{MX}-\mathbf{YP} \rVert_{2}^{2} 
\end{equation}

When neither the correspondences, nor the linear transformation are known the problem takes the form as in Eqn. (\ref{psr})

\begin{equation}
\label{psr}
\underset{\mathbf{M}\in O(d),\mathbf{P}\in\mathbb{\mathcal{P}}_{n}}{\arg\min}\lVert \mathbf{MX}-\mathbf{YP} \rVert_{2}^{2}
\end{equation}

There have been attempts where graph matching algorithms (the approximate form of graph isomorphism) are used to solve point set registration optimization problems~\cite{4587538}. However, it is not very common to observe methods going the other way around. Namely, casting a graph matching problem as a point set registration problem is not a conventional idea. In this paper, we show that graph isomorphism can be cast as an exact rigid point set registration problem in $\mid V \mid$ dimensions using a simplex embedding and sampling. Therefore, existence of an exact solution to such an instance is equivalent of solving the original graph isomorphism decision problem.

The rest of this paper is organized as follows. Section~\ref{form} introduces the main theory of casting a graph isomorphism instance as a point set registration problem. Section ~\ref{vars} then investigates the introduced theory in terms of certain variations such as in the cases of automorphisms, weighted graphs, or hypergraphs to name a few. The rest of the formulation includes using the optimization equation introduced as a distance measure between two graphs as in Section~\ref{dist}. Then, Section~\ref{related} examines the related literature and tries to make connections between our formulation and related concepts. Finally, Section~\ref{discuss} presents a discussion and a conclusion of this study.

\section{Formulation}
\label{form}
A graph is an abstract data type consisting of two sets, namely a set of vertices $V$, and a set of edges $E$ (i.e. unordered pairs of the vertices). Until further notice, by a graph a simple graph is meant, namely the graph is undirected, unweighted, and does not contain any loops or multiple edges. 

Conventionally, graphs are depicted on planes, where vertices are represented by disks and edges are represented by continuous lines linking a pair of disks (vertices). Although such representation is needed for graphs to be drawn on paper, such a depiction limits further processing of the graph. It is possible to come up with more expressive geometric representations of graphs in higher dimensions.

\subsection{Simplex embedding}

The intuition behind a meaningful depiction of a graph is the fact that $V$ is a set and there is no ranking imposed. Therefore, in a meaningful depiction vertices must be equidistant from each other. The first nontrivial case is when $\mid V \mid =3$ and in that case vertices must coincide with the vertices of an equilateral triangle to be equidistant from each other. When $\mid V \mid =4$, then a planar depiction is not possible and we need vertices to coincide with the vertices of a regular tetrahedron. With each new vertex, we need to increase the dimension by one to keep the vertices equidistant from each other, arriving at a formulation such as a regular simplex.
In this manner, we can embed the set $V$ onto the vertices of the standard $\mid V \mid -1$-simplex in $\mathbb{R}^{\mid V \mid}$ where $\mathbf{S}_{V} = \mathbf{I}$ corresponds to the coordinates of the new vertices in this embedding, the identity matrix having size $\mid V \mid \times \mid V \mid$. Since the standard simplex is not centered on the origin, additional step is needed to center it and get rid of any translation if we are to use orthogonal-only transformations in our formulation. This is achieved by calculating the mean of points and subtracting it from each column as in Alg.\ref{simplex} written in Matlab pseudo-script style.

\begin{algorithm}[t]
\footnotesize
    \caption{\small A function to acquire vertices of a centered standard simplex representing the set of vertices V}
		\label{simplex}
    \begin{algorithmic}[1]
    \Function{SimplexPoints}{$V$}
		\State $d \gets \mid V \mid$
		\State $S_{V} \gets eye(d,d)$
		\State $S_{V} \gets S_{V}-repmat(mean(S_{V},2),1,d)$
		\State \Return $S_{V}$
    \EndFunction
    \end{algorithmic}
    \end{algorithm} 

\subsection{Sampling}

The edges of the graph are still assumed to be continuous entities and trace the edges of the underlying simplex. However, in order to convert the graph into a point cloud, edges must be represented as points also. Therefore, a sampling is possible where the midpoint of two vertices is included as an additional point depending on the existence of an edge in-between. Instead of a single sample point, multiple sample points can also be attained. However, a single sample point is enough to assure the existence of an edge. This sampling procedure is given in Alg.\ref{edges}.

\begin{algorithm}[t]
\footnotesize
    \caption{\small A function to acquire sample points representing the edges of the graph}
		\label{edges}
    \begin{algorithmic}[1]
    \Function{EdgePoints}{$S_{V},E$}
		\State $S_{E} \gets zeros(size(S_{V},1),size(E,1))$
		\For{$i \gets 1$ to $size(E,1)$} 
				\State $e \gets E(i,:)$
        \State $S_{E}(:,i) \gets (S_{V}(:,e(1))+S_{V}(:,e(2)))/2$
		\EndFor
		\State \Return $S_{E}$
    \EndFunction
    \end{algorithmic}
    \end{algorithm}

\subsection{Casting as a point set registration problem}

In an isomorphism test, we have a pair of graphs that can be denoted as $G_{1}$ and $G_{2}$, where $V_{1}$ and $E_{1}$ denote the vertices and edges of the graph $G_{1}$, and $V_{2}$ and $E_{2}$ denote the vertices and edges of graph $G_{2}$. It is now possible for us to convert $G_{1}$ into a point cloud form by attaining $\mathbf{S}_{V_{1}}$ and $\mathbf{S}_{E_{1}}$ and concatenating them by $\mathbf{S}_{1} = \left[\mathbf{S}_{V_{1}} \: \mathbf{S}_{E_{1}}\right]$ to obtain $\mathbf{S}_{1}$ where it is a $\mid V_{1} \mid \times \left[\mid V_{1} \mid +\mid E_{1}\mid \right]$ matrix.
Without explicitly formulating, we can then denote the second graph in its point cloud form as $\mathbf{S}_{2} = \left[\mathbf{S}_{V_{2}} \: \mathbf{S}_{E_{2}}\right]$ having a size of $\mid V_{2} \mid \times \left[\mid V_{2} \mid + \mid E_{2} \mid\right]$
Having two point clouds we can then use Eqn.(\ref{psr}) to formulate graph isomorphism of $G_{1}$ and $G_{2}$ in their point set forms of $\mathbf{S}_{1}$ and $\mathbf{S}_{2}$ respectively as in Eqn.(\ref{iso}). Note that, an exact registration is needed for two graphs to be isomorphic.  

\begin{equation}
\label{iso}
\lVert \mathbf{MS_{1}}-\mathbf{S_{2}P} \rVert_{2}^{2} = 0 \quad s.t. \quad \mathbf{M}\in O(d) \wedge \mathbf{P}\in\mathbb{\mathcal{P}}_{n}
\end{equation}

\subsection{Special consideration of the transformation matrix}

Note that, the required transformation matrix $\mathbf{M}$ covers all the possible orthogonal real matrices, thus there are infinitely many candidates. A careful consideration suggests that the set of valid transformations is finitely many. Note that, the overall shape of the point cloud entity is a centered standard simplex for both $\mathbf{S}_{1}$ and $\mathbf{S}_{2}$. Therefore, a valid transformation should map the vertices of the simplex onto the same shape again. In other words, the symmetry group of a regular $(\mid V_{1} \mid -1)$-simplex is the symmetric group $\mathbf{Sym}(\bf{S}_{V_{1}})$. Therefore, there are only a fixed number of valid configurations.

\subsection{Special consideration of the permutation matrix}

A special consideration of the permutation matrix $\mathbf{P}$ is also needed. Note that, $\mathbf{S}_{1}$ is a concatenation of $\mathbf{S}_{V_{1}}$ and $\mathbf{S}_{E_{1}}$, similar structure is observed for $\mathbf{S}_{2}$. Therefore, $\mathbf{P}$ must be designed in such a way that vertices ($\mathbf{S}_{V_{1}}$) are mapped to vertices ($\mathbf{S}_{V_{2}}$), and vertices corresponding to edges ($\mathbf{S}_{E_{1}}$) must be mapped to vertices corresponding to edges($\mathbf{S}_{E_{2}}$). A more careful consideration suggests that one may not solely change the mapping of edge vertices, but can only alter the mapping of vertices. Given an order of vertices, the order of edge vertices should be automatically set. This poses a special restriction to the general setting of $\mathbf{P}\in\mathbb{\mathcal{P}}_{n}$, and can be overwritten with $\mathbf{P}\in\mathbb{\mathcal{P}}_{n}^{VE}$  where $\mathbb{\mathcal{P}}_{n}^{VE}$ designates the restricted form of permutation matrix set. With these considerations we arrive at the final formulation as in Eqn.(\ref{isofin})

\begin{equation}
\label{isofin}
\lVert \mathbf{MS_{1}}-\mathbf{S_{2}P} \rVert_{2}^{2} = 0 \quad s.t. \quad \mathbf{M}\in \mathbf{Sym}(\bf{S}_{V_{1}}) \wedge \mathbf{P}\in\mathbb{\mathcal{P}}_{n}^{VE}
\end{equation}

\section{Variations}
\label{vars}
To gain further insight and measure the capabilities of the proposed formulation, certain variations of the concept are to be investigated.

\subsection{Automorphisms}
An automorphism of a graph $G$ is a graph isomorphism from $G$ to itself. The graph automorphism problem (GA) tests whether a graph has a nontrivial automorphism (i.e. excluding identity). It is in class NP. However, similar to GI, it is still unknown whether it has a general polynomial time algorithm or it is NP-complete. It is also known that GA is polynomial-time many-one reducible to GI, but the reduction on the opposite way is unknown~\cite{toran2004hardness}.

There is an intricate relation between GA and GI. For example, two graphs $G_{1}$ and $G_{2}$ are isomorphic if and only if the disjoint union of them has an automorphism that sends these two components to each other~\cite{LUKS198242}. 

At this point, let us present GA in point set form using our formulation. Since there is one graph as input there will also be a single point set $\mathbf{S}$. Therefore, Eqn.(\ref{isofin}) simplifies to Eqn.(\ref{auto2}). 


\begin{equation}
\label{auto2}
\begin{split}
\mathbf{MSP^{T}} = \mathbf{S} \quad s.t. \quad \mathbf{M} &\in \mathbf{Sym}(\bf{S}_{V}) \setminus \mathbf{I} \quad \wedge \\
\mathbf{P}&\in\mathbb{\mathcal{P}}_{n}^{VE} \setminus \mathbf{I}
\end{split}
\end{equation}

Note that, $\mathbf{MSP^{T}}=\mathbf{S}$ can be interpreted both as $\mathbf{(MS)P^{T}}$ or $\mathbf{M(SP^{T})}$ due to associativity of matrix multiplication. Semantically, it means that after a change in correspondence of the points with $\mathbf{P^{T}}$, we can transform the point set into its original form through the geometric transformation $\mathbf{M}$. In simpler terms, this is a representation of the symmetry group of the point set $\mathbf{S}$. If the order (or cardinality) of this group is more than 1, that corresponds to $G$ having at least one nontrivial automorphism. This question then boils down to finding whether the order of the symmetry group of a set of $\mid V \mid + \mid E \mid$ points in $\mid V \mid$ dimensions is bigger than one. Counting the number of automorphisms (\#GA) is then equivalent to finding the order of this symmetry group. 

A special note is that this symmetry group must be a subgroup of the all possible mappings of a $(\mid V \mid -1)$-dimensional simplex onto itself (i.e. the symmetric group of order $\mid V \mid$). For example, in the case of the input graph being a complete one, the symmetry group and the symmetric group coincide.

\subsection{Graph generalizations}

It is possible to generalize the concept of a simple graph in different directions. In a directed graph, a set of ordered pairs of vertices now describes the set of edges. By sampling the further midpoint of the original midpoint and the second vertex consistently, it is possible to distinguish the direction of edges in the point set form. A hypergraph is a different generalization of the graph concept, where an edge can now relate an arbitrary number of vertices. It is already shown that, GI is polynomial-time equivalent to HI~\cite{Arvind2015}. In our context, the vertices of the hypergraph again take the places of vertices of a simplex. Then, each hyperedge can be converted to a point as the mean point of the vertices it relates. However, for weighted graphs our sampling process is possibly not suitable. For example, in the case of real weights or negative weights it is hard to come up with a meaningful sampling process.

\subsection{Subgraph isomorphism}
\label{sgi}
In subgraph isomorphism problem, it is checked whether a graph $G_{1}$ contains a subgraph that is isomorphic to the other graph $G_{2}$. This variant is of special interest as subgraph isomorphism problem is known to be NP-complete~\cite{10.1145/800157.805047}. Therefore, the point set form of it is also expected to be more involved. First of all, the number of vertices of $G_{1}$ must be at least that of $G_{2}$. Therefore, it is better to embed $G_{1}$ onto a simplex first and determine the dimensionality of the space. $G_{1}$ can be conventionally converted to its point cloud form $\mathbf{S}_{1}$. Similarly, $G_{2}$ can then be processed as points on a simplex of dimensions $\mid V_{2} \mid -1$ in a space of $ \mid V_{1} \mid$ dimensions, where both point clouds are centered on the origin. 

A striking observation is that orthogonal-only transformations will not be enough to find a perfect fitting of $\mathbf{S}_{2}$ onto $\mathbf{S}_{1}$ as possible substructures of $S_{1}$ are not centered around the origin but has translational offsets. Note that, homogeneous coordinates provide a convenient way to represent a translation as a matrix multiplication~\cite{Li2001}. Therefore, in our extended formulation, we denote such a transformation as $\mathbf{D}$. 

There is another problem needed to be resolved for a matrix representation. Although the dimensions of $\mathbf{S}_{1}$ and $\mathbf{S}_{2}$ match, namely their row sizes are equal, their column sizes or the number of points they have may not be equal. Therefore, a combination of proper size is needed to select the columns of $\mathbf{S}_{1}$, denoted as a list of indexes $\mathbf{c}$. Note that, the indexes can be in any order, since if they are correct there is also the permutation matrix $\mathbf{P}$ to find the correct correspondence. Therefore, the sole function of this indexing is to reduce the size of point count of $\mathbf{S}_{1}$ to that of $\mathbf{S}_2$. Note that, $\mathbf{c}$ is defined in such a way that vertices are selected from vertices, and edge vertices are selected from edge vertices. With these issues addressed, subgraph isomorphism problem can be represented in point cloud form as in Eqn.(\ref{subg})

\begin{equation}
\label{subg}
\begin{split}
&\lVert \mathbf{DMS_{1}\big[:,c\big]}-\mathbf{S_{2}P} \rVert_{2}^{2} = 0 \quad s.t. \\
&\mathbf{M} \in O(d) \\
&\mathbf{P}\in\mathbb{\mathcal{P}}_{n}^{VE}  \\
&\mathbf{D} \in \mathbb{D} \\
&\mathbf{c} \in {V_{1} \choose \mid V_{2} \mid} \cup {E_{1} \choose \mid E_{2} \mid } \\
\end{split}
\end{equation}

An important note is that the transformation matrix $\mathbf{M}$ has finitely many configurations, but it is hard to formulate it so it is reverted back to the more general setting of orthogonal group $O(d)$.

\subsection{Alternative vector set representation}
Throughout this paper, we used point sets. Points can also be thought of as vectors. Namely, instead of aligning points in the least squares sense, an alignment of vectors can alternatively be formulated. There are two possible repercussions in that case. Firstly, normalizations will most possibly be needed. Secondly, such an alternative will not be able to formulate subgraph isomorphism problem as in Section~\ref{sgi}, because it will not be possible to model additional translational transformation. Therefore, point cloud formulation seems as a more general representation for graphs in general.

\section{A distance notion of two graphs}
\label{dist}
In the case of non-isomorphism, the optimization version of the problem, namely Eqn.(\ref{iso}) in the form of Eqn.(\ref{psr}), can be used to introduce a notion of distance between two graphs having the same number of vertices and edges. We call this notion as Graph Geometric Distance (GGD). If this distance is zero, it means that two given graphs are isomorphic. It may be important to investigate the other extreme when one of the graphs remains fixed, namely the maximum possible distance to a graph among all candidates. Such a notion can be grasped with the term telomorphism. Namely, a graph (fixed) is telomorphic to another if the distance between them is the maximum possible among all other candidates. It is important to note that a telomorphism can also be defined over a different distance notion. Such observation is important if one wants to investigate the relationship between different distance notions of two graphs.  

The idea of telomorphism is closely related to equivalence classes that graph isomorphism relation imposes. However, without a distance notion it is not possible to comment on the relationship between these equivalence classes. Studying relationships between equivalence classes may be important in the quest for graph canonization, a tightly related concept. A solution to graph canonization problem implies a solution to the graph isomorphism problem as one can first compute the canonical forms of two graphs and then check whether these two forms are identical. Formally, the graph canonization problem is at least as computationaly hard as the graph isomorphism problem, but it is still not known whether these two problems are polynomial time equivalent~\cite{10.1007/978-3-540-77120-3_71}. Optimistically, an efficient geometric canonical form of simple graphs may be deduced from the representation in this study.

\section{Related work}
\label{related}

\subsection{On graph-simplex correspondence}

There exists a different graph-simplex correspondence where every connected, undirected, unweighted graph of $\mid V \mid$ nodes corresponds to one specific simplex in $\mid V \mid -1$ dimension\cite{fiedler1993geometric,10.1093/comnet/cny036}. In that formulation, each graph is mapped to a whole simplex with varying edge-sizes depending on their vertex similarities~\cite{10.1093/comnet/cny036}. It is questionable whether such a formulation is more suitable for GI. For example, it is not clear whether two graphs are isomorphic if and only if they have the same simplex structure. Nevertheless, such a correspondence is powerful enough to find connections between certain graph characteristics and simplex properties, but enforcing connectedness is a limitation.

\subsection{Application areas}

To gain further insight one may investigate the related application domain of a certain problem. Traditionally, graph isomorphism is used as the basis of Layout versus Schematic (LVS) circuit design step~\cite{10.5555/800261.809095}. In a more modern context, graph isomorphism in its original form is primarily applied in cheminformatics to identify a chemical compound within a chemical database. Namely, it is important to designate a canonical representation for each chemical compound. In this line, The Simplified Molecular-Input Line-Entry System (SMILES) and International Chemical Identifier (InChl) have been designed to provide standard ways of encoding molecular information based on the canonization of the graph representing the molecule~\cite{10.1145/2366316.2366334}. 

On the other hand, point set registration in perceivable dimensions has applications in computer vision~\cite{7299195} and robotics~\cite{7747236}. Recently, higher dimensional applications also include natural language processing, where word embeddings are used~\cite{pmlr-v89-grave19a}. Therefore, our formulation is important in that it relates the first domain of applications with the second one. In this regard, an efficient solution to graph isomorphism or point set registration would mean an advance in many related fields.

\subsection{On distance notion of two graphs}
\label{ddist}
Although a distance within a graph is a well-studied concept, the distance notion between two graphs remains elusive. By distance between two graphs usually Graph Edit Distance (GED) is meant~\cite{6313167}. GED algorithms are categorized based on the graphs being attributed or not and also on the definition of costs for edit operations~\cite{Gao2010}. There is no trivial relationship between GED and GGD introduced in Section ~\ref{dist}. First, a GGD must be introduced for two arbitrary graphs instead of graphs having the same number of vertices and edges. In that case, the optimization form of Eqn.~(\ref{subg}) (subgraph isomorphism) can be adopted. 
However, since subgraph isomorphism is NP-complete, such a distance notion may not be practical, but may still have theoretical significance. 

A more systematic review of various distance notions for graph comparison is given by~\cite{10.1371/journal.pone.0228728}. It is noted that, the distances must scale linearly or near-linearly to be practical in real world scenarios. From a practical perspective distance notions are studied under 5 headings in~\cite{10.1371/journal.pone.0228728}. These include spectral, matrix, and feature-based distances as 3 conventional ways. Learning based approaches exist where algorithms learn an embedding from a set of graphs into the Euclidean space and then compute a notion of similarity in between~\cite{GOYAL201878,pmlr-v97-li19d}. Most frequently, graph neural networks~\cite{4700287} are utilized in learning based approaches. In a recent study, a neural embedding framework named graph2vec is used to learn vector representations of arbitrary sized graphs~\cite{narayanan2017graph2vec}. However, such study does not refer to graph isomorphism problem at all. Namely, it must be the case that two isomorphic graphs give out the same vector representation. Within the machine learning perspective, a successful neural network must output the same results for any pair of isomorphic graphs. In extended words, graph isomorphism problem can be cast as an ultimate machine learning clustering problem where given a graph of specified size the machine should be able to determine the correct cluster (or equivalence class) of the graph. In our point set formulation of a graph, such a machine must then be invariant to relevant transformations in $\mid V \mid$ dimensions. 

\section{Discussion and Conclusion}
\label{discuss}

A more traditional way of representing a graph is through an adjacency matrix. Using this formulation, it is possible to show that two graphs $G_{1}$ and $G_{2}$ are isomorphic if and only if $\bf{A_{1}} = \bf{P^{T}}\bf{A_{2}}\bf{P}$ holds where $\bf{A_{1}}$ and $\bf{A_{2}}$ are adjacency matrix representations of two graphs and $\bf{P}$ is a $\mid V \mid \times \mid V \mid$ permutation matrix. This appears to be a simpler formulation than Eqn.(\ref{isofin}) for graph isomorphism. Namely, a single permutation matrix is enough as a solution. However, in the literature there happens to be not much efficient progress towards graph isomorphism using the adjacency matrix form. 

Even though a single permutation matrix is enough in adjacency matrix form, how is it possible that a point set formulation requires both a transformation matrix and a permutation matrix as in Eqn.(\ref{isofin})? Actually, in our formulation, a permutation matrix is only needed to give out the exact correspondence between the nodes of the graphs if there exists an isomorphism. However, to tell whether there exists an isomorphism, after the required transformation, checking whether two point sets have the same positions included is enough. If such a condition arises, it is then trivial to find the permutation matrix. Therefore, decision version of the problem and finding a correspondence if it exits are polynomially equivalent~\cite{MATHON1979131}. 

Moreover, it is possible to refer to Cartan--Dieudonn{\'e} theorem~\cite{Gallier2011} for the construction of the required transformation matrix. In other words, it is possible to construct any candidate transformation through a series of reflections. In simple terms, we can swap the place of two simplex vertices by reflecting the simplex through the flat passing through the remaining vertices and the midpoint between these two simplex vertices. With this single primitive operation, it is possible to represent any needed transformation. 

Another note is related to the subgraph isomorphism problem. Most of these distance notions mentioned in this study do not satisfy the conditions of a metric. For example, in a related manner, the distance notion for comparing two different sized graphs introduced in Section~\ref{ddist} does not satisfy the triangle inequality. It seems hard to come up with a distance notion of two arbitrary graphs that is also a metric. For example, based on the subgraph isomorphism problem, the distance of whole to a part may be designated as zero. Then, the distance to another differently shaped part may then be zero. However, the distance between these two parts must be zero for triangle inequality to hold, but these two parts are not isomorphic. Therefore, this opens up a further discussion on the distance notion of a part-whole hierarchy.  

It is important to repeat the observation that in the general setting namely in comparison of two differently sized graphs such as in the case of subgraph isomorphism, vector set representation might be limited as translations are also needed to find a correspondence. Therefore, it seems that representing a graph in a point set form is a more general approach. A related future work includes determining the exact finite subset of orthogonal transformations needed for subgraph isomorphism problem given two candidate graphs.  

As a final note, this representation might pave way to an efficient graph canonization. In this regard, graph invariants play a crucial role. A graph invariant is complete if the identity of invariant forms of any two graphs implies their isomorphism. For example, Zernike moments can be utilized in 2D for rotation-invariant feature extraction of images~\cite{55109}. The problem then is to generalize this concept to efficient high-dimensional orthogonal transformation invariant feature extraction on the point set form of graphs. Namely, this point set representation can be an additional layer before casting a graph as a vector through feature extraction. In this light, complete graph invariants based on this point set form might be important both theoretically and practically.

\bibliography{mybibfile}

\end{document}